\documentclass[10pt]{article}
\title{{ On some expansions for
the Euler Gamma function and the Riemann Zeta
function. }}
\author{Grzegorz Rz\c{a}dkowski
}
\usepackage{amssymb}
\usepackage{latexsym}
\usepackage{amsthm}
\usepackage{graphics}
\date{Faculty of Mathematics and Natural Sciences,\\ Cardinal Stefan Wyszy\'nski 
University in Warsaw,\\ Dewajtis 5, 01 - 815 Warsaw, Poland\\ g.rzadkowski@uksw.edu.pl  \\
grzerzad@gmail.com} 
\begin{document}
\maketitle
\newtheorem{lemma}{Lemma}
\newtheorem{theorem}{Theorem}
\newtheorem{remark}{Remark}
\newtheorem{statement}{Statement}
\def \bangle{ \atopwithdelims \langle \rangle}
\textbf{The paper has been published in:\\
J. Comp. Appl. Math. 236 (2012) pp. 3710-3719.}
\begin{abstract}
In the present paper we introduce some expansions, based on the falling factorials, for the Euler Gamma function and the Riemann Zeta function. In the proofs we use the Fa\'a di Bruno formula, Bell polynomials, potential polynomials, Mittag-Leffler polynomials, 
derivative polynomials and special numbers (Eulerian numbers and Stirling
numbers of both kinds). We investigate the rate of convergence of the series and give some numerical examples.
\end{abstract}
\noindent 2010 {\it Mathematics Subject Classification}:
11B83, 11M06 .\\
\noindent \emph{Keywords: } Euler Gamma function, Riemann Zeta function, Bell polynomials, potential polynomials,  Mittag-Leffler polynomials, derivative polynomials.
\section{Introduction}
Let us first recall some basic facts concerning special numbers and expansions. By $\displaystyle {n \brack k}$ we denote the Stirling number of the first kind (number of the ways of partitioning a set of $n$ 
elements into $k$ nonempty cycles, see \cite{GKP}). It is set $\displaystyle { n \brack 0} =0 $ if $n>0$, $\displaystyle { 0 \brack 0} =1 $,  $\displaystyle { n \brack k}=0 $ for $k>n$ or $k<0$. The Stirling numbers of the first kind fulfill the recurrence formula
\begin{equation}\label{re}
	{ n \brack k} =(n-1){ n-1 \brack k} + { n-1 \brack k-1}.
\end{equation}
If $n>k$ then using formula (\ref{re}) in each step to the last term of the resulting sum we get
\begin{equation}\label{re2}
	{ n \brack k} = \sum\limits_{j=1}^{k}(n-j){ n-j \brack k+1-j}.
\end{equation}
Stirling numbers of the first kind have the following generating function (see Comtet \cite{C}, p. 50, p. 135)
\begin{equation}\label{gf}
	(1-t)^{-u}=1+\sum\limits_{1\le k\le n}{ n \brack k}\;\frac{t^{n}}{n!}\;u^{k}.
\end{equation}
We use common notations for the falling factorial 
	\[(x)_{k}=x(x-1)\cdots (x-k+1)
\]
 and for the rising factorial (Pochhammer's symbol)
	\[x^{(k)}=x(x+1)\cdots (x+k-1).
\]
We denote by $B_{n,k}=B_{n,k}(x_{1},x_{2},\ldots, x_{n-k+1})$ (see \cite{B}, \cite{N}, \cite{C} (p.133)) the exponential partial Bell polynomials in infinite number of variables $x_{1},x_{2},x_{3},\ldots$. The polynomials are defined by the formal double series expansion in variables $t$ and $u$
\begin{equation}\label{Bp}
	\exp\biggl(u\sum\limits_{m\ge 1}x_{m}\frac{t^{m}}{m!}  \biggr)=
	1+\sum\limits_{n\ge 1}\frac{t^{n}}{n!}\biggl\{\sum\limits_{k=1}^{n}u^{k}B_{n,k}(x_{1},x_{2},\ldots )   \biggr\}. 
\end{equation}
We denote by $P_{n}^{r}$ the potential polynomials (see \cite{C}, Theorem B p. 141) which are defined for each complex number $r$ by
\begin{equation}\label{pp}
	\biggl(1+\sum\limits_{n\ge 1}g_{n}\frac{t^{n}}{n!}\biggr)^{r}=1+\sum\limits_{n\ge 1}P_{n}^{r}\frac{t^{n}}{n!}
\end{equation}
and 
\begin{equation}\label{pp2}
	P_{n}^{r}=P_{n}^{r}(g_{1},g_{2},\ldots,g_{n})=\sum\limits_{1\le k \le n}(r)_{k}B_{n,k}(g_{1},g_{2},\ldots ). 
\end{equation}
Formula (\ref{pp}) is a particular case of the Fa\'a di Bruno formula and $P_{n}^{r}$ (given by (\ref{pp2})) is the $n$th derivative (in a point $x=a$) of the function $(G(x))^{r}$, where $G(x)$ is given as the convergent power series $G(x)=1+\sum\limits_{n\ge 1}g_{n}t^{n}/n!, \quad t=x-a,\; G(a)=1$.\\
We investigate expansions, which involve the falling factorials, for the Euler Gamma function and for the integral (\ref{Z}). The last is expressed in terms of the Riemann Zeta function. For the coefficients of our series we give simple recurrence formulae. Some series for the the Riemann Zeta function based on falling factorials have been studied, for example, by Flajolet and Vepstas \cite{FV}, who demonstrated the importance of such expansions. The coefficients of their expansion are expressed in terms of the values of the Zeta function in integers. \\
The article is organized as follows. In Section 2 we present the construction of the expansion, based on falling factorials, for the Euler Gamma function. For its coefficients we give the recurrence formula and an explicit formula involving the Stirling numbers of the first kind. In Section 3 we present the basic properties of the derivative polynomials which are used in the next sections. Section 4 is devoted to the construction of the expansion for integral (\ref{Z}). For the coefficients of the expansion we give the recurrence formula and an explicit formula, which uses the coefficients of the Mittag-Leffler polynomials. In Section 5 we examine the rate of convergence of the series introduced in Sections 2 and 4. Moreover, we show the results of two numerical experiments. The paper is concluded in Section 6.   
\section{The Euler Gamma function}
Substituting in the integral
	\[\Gamma (s+1)=\int\limits_{0}^{\infty}x^{s}e^{-x}dx
\]
$x=-\log(1-t)$ we have 
\begin{equation}\label{G}
	\Gamma (s+1)\!=\!\!\int\limits_{0}^{1}\!(-\log(1-t))^{s}\!dt\!=\!\!\int\limits_{0}^{1}\!t^{s}\left(\frac{1}{t}\log\frac{1}{1\!-\!t}\right)^{s}\!\!dt\!=\!\!\int\limits_{0}^{1}\!t^{s}\left(1\!+\!\frac{t}{2}\!+\!\cdots\right)^{s}\!\!dt.
\end{equation}
Our first aim is to find the values of the Bell polynomials $B_{n,k}$ for the sequence $(1/2,2!/3,3!/4,\ldots)$. Using expansion (\ref{Bp}) we get
\begin{eqnarray*}
\lefteqn{\exp\biggl\{u\left(\frac{1}{2}t+\frac{1}{3}t^{2}+\frac{1}{4}t^{3}+\cdots\right)\biggr\}}\\
&&=e^{-u}\exp\biggl\{u\left(1+\frac{1}{2}t+\frac{1}{3}t^{2}+\frac{1}{4}t^{3}+\cdots\right)\biggr\}\\
&&=e^{-u}\exp\biggl\{\frac{u}{t}(-\log(1-t))\biggr\}=e^{-u}(1-t)^{-\frac{u}{t}}\\
&&=e^{-u}\biggl\{1+\sum\limits_{1\le k\le n}{n \brack k}\frac{t^{n}}{n!}\left(\frac{u}{t}\right)^{k}\biggr\}
= e^{-u}\biggl\{1+\sum\limits_{1\le k\le n}{n \brack k}\frac{t^{n-k}u^{k}}{n!}\biggr\}\\
&&=\left(1-\frac{u}{1!}+\frac{u^{2}}{2!}-\frac{u^{3}}{3!}+\frac{u^{4}}{4!}-\cdots\right)+
\sum\limits_{j=0}^{\infty}\sum\limits_{1\le k\le n}(-1)^{j}{n \brack k}\frac{t^{n-k}u^{k+j}}{n!j!}.
\end{eqnarray*}
Putting $n-k=\alpha \ge 0$, $k+j=\beta\ge 1$, ($n=k+\alpha,\;j=\beta-k$) we see that the coefficient of $t^{\alpha}u^{\beta}$ is
\begin{equation}
	\sum\limits_{k=1}^{\beta}{k+\alpha \brack k}\frac{(-1)^{\beta-k}}{(k+\alpha)!(\beta-k)!}=
	\frac{1}{(\alpha+\beta)!}\sum\limits_{k=1}^{\beta}(-1)^{\beta-k}{k+\alpha \brack k}{\alpha+\beta \choose \beta-k}.
\end{equation}
By denoting
	\[c_{\alpha, \beta}=\sum\limits_{k=1}^{\beta}(-1)^{\beta-k}{\alpha+k \brack k}{\alpha+\beta \choose \beta-k},
\]
we write the value of the Bell polynomial $B_{\alpha,\beta}$ for the sequence \\ $(1/2,2!/3,3!/4,\ldots)$ as 
\begin{equation}\label{Bp2}
	B_{\alpha,\beta}=\frac{\alpha!}{(\alpha+\beta)!}\;c_{\alpha,\beta}.
\end{equation}
The following Table \ref{tab:tb1} gives the first few values of
$c_{\alpha,\beta}$.\\
\setlength{\belowcaptionskip}{4pt}
\begin{table}[h]
\caption{Coefficients $c_{\alpha,\beta}$}
\centering
$        \begin{array}{cccccccc}
        \alpha \backslash \beta &\  0 &\  1 & 2 & 3 & 4 & 5 &\  6    \\
        0 &\  1 &\  0 & &&&& \\
        1 &\  0 &\  1 & 0  &&&& \\
        2 &\  0 &\  2 & 3 &0&&& \\
        3 &\  0 &\  6 & 20 & 15 &0&& \\
        4 &\  0 &\  24 & 130 & 210 & 105 &0&\\
        5 &\  0 &\  120 & 924 & 2380 & 2520 & 945 &\ 0\\
\end{array}
$
\label{tab:tb1}

\end{table}
\begin{lemma}
Numbers $c_{\alpha,\beta}$ fulfill the recurrence formula
\begin{equation}\label{cnk}
	c_{\alpha,\beta}=(\alpha+\beta-1)(c_{\alpha-1,\beta}+c_{\alpha-1,\beta-1}).
\end{equation}
\end{lemma}
\begin{proof}
Putting into (\ref{cnk}) 
{\setlength\arraycolsep{2pt}
\begin{eqnarray*}
c_{\alpha-1, \beta}&=&\sum\limits_{k=1}^{\beta}(-1)^{\beta-k}{\alpha+k-1 \brack k}{\alpha+\beta-1 \choose \beta-k},\\
c_{\alpha-1, \beta-1}&=&\sum\limits_{k=1}^{\beta-1}(-1)^{\beta-1-k}{\alpha+k-1 \brack k}{\alpha+\beta-2 \choose \beta-1-k},
\end{eqnarray*}}
and adding similar terms we see that formula (\ref{cnk}) is equivalent to
\begin{equation}\label{cnka}
	c_{\alpha,\beta}= (\alpha+\beta-1)\sum\limits_{k=1}^{\beta}(-1)^{\beta-k}{\alpha+k-1 \brack k}{\alpha+\beta-2 \choose \beta-k}.
\end{equation}
From the other side 
{\setlength\arraycolsep{2pt}
\begin{eqnarray}
 c_{\alpha, \beta}&=&\sum\limits_{k=1}^{\beta}(-1)^{\beta-k}{\alpha+k \brack k}{\alpha+\beta \choose \beta-k}\nonumber\\
 &=&
 \sum\limits_{k=1}^{\beta}(-1)^{\beta-k}{\alpha+\beta \choose \beta-k}
 \sum\limits_{i=1}^{k}(\alpha+k-i){\alpha+k-i \brack k+1-i}\nonumber\\
 &=&\sum\limits_{m=1}^{\beta}(\alpha+\beta-m){\alpha+\beta-m \brack \beta+1-m}\sum\limits_{i=1}^{m}
 (-1)^{m-i}{\alpha+\beta \choose m-i}\nonumber\\
 &=& \sum\limits_{m=1}^{\beta}(-1)^{m-1}(\alpha+\beta-m){\alpha+\beta-1 \choose m-1}{\alpha+\beta-m \brack \beta+1-m}\nonumber\\
 &=& (\alpha+\beta-1)\sum\limits_{m=1}^{\beta}(-1)^{m-1}{\alpha+\beta-2 \choose m-1}{\alpha+\beta-m \brack \beta+1-m},\label{lst}
\end{eqnarray}}
where we introduced the new parameter $m$ by formula $m-i=\beta-k$. If we put $\alpha +k-1=\alpha+\beta-m$ (i.e. $k-1=\beta-m$) 
the formula (\ref{lst}) converts to (\ref{cnka}).\\
At the beginning of the above calculation we used the formula (\ref{re2}) and then formulae
{\setlength\arraycolsep{2pt}
\begin{eqnarray*}
&&\sum\limits_{i=1}^{m}
 (-1)^{m-i}{n \choose m-i}=(-1)^{m-1}{n-1 \choose m-1},\\
&&(n-m){n-1 \choose m-1}=(n-1){n-2 \choose m-1},
\end{eqnarray*}}
which are easy to verify.
\end{proof}
Therefore using (\ref{pp}) and (\ref{Bp2}) we write the integrand in (\ref{G}) in the form
\begin{eqnarray}
\lefteqn{(-\log(1-t))^{s}=t^{s}\left\{1+t\frac{s}{2}+t^{2}\left(\frac{2}{3!}s+\frac{3}{4!}s(s-1)\right)\right.
}\nonumber\\
&&\left.+
t^{3}\left(\frac{6}{4!}s+\frac{20}{5!}s(s-1)+\frac{15}{6!}s(s-1)(s-2)\right)+\cdots\right\}\nonumber\\
&&=t^{s}\left\{1+\sum\limits_{\alpha=1}^{\infty}\frac{t^{\alpha }}{\alpha!}\sum\limits_{\beta=1}^{\alpha}(s)_{\beta}B_{\alpha,\beta}  \right\}=
t^{s}+\sum\limits_{\alpha=1}^{\infty}t^{\alpha +s}\sum\limits_{\beta=1}^{\alpha}(s)_{\beta}\frac{c_{\alpha,\beta}}{(\alpha+\beta)!}\label{Ga}
\end{eqnarray}
and finally for the integral (\ref{G}) we get the following expansion 
\begin{eqnarray}
\lefteqn{\Gamma (s+1)=\frac{1}{s+1}+\frac{1}{s+2}\cdot \frac{s}{2}+\frac{1}{s+3}\left(\frac{2}{3!}s+\frac{3}{4!}s(s-1)\right)
}\nonumber\\
&&+\frac{1}{s+4}\left(\frac{6}{4!}s+\frac{20}{5!}s(s-1)+\frac{15}{6!}s(s-1)(s-2)\right)+\cdots \nonumber\\
&&=\frac{1}{s+1}+\sum\limits_{\alpha=1}^{\infty}\frac{1}{\alpha +s+1}\sum\limits_{\beta=1}^{\alpha}(s)_{\beta}\frac{c_{\alpha,\beta}}{(\alpha+\beta)!}.\label{Gam}
\end{eqnarray}
By denoting
	\[g_{\alpha,\beta}=g_{\alpha,\beta}(s)=(s)_{\beta}\frac{c_{\alpha,\beta}}{(\alpha+\beta)!},
\]
and using (\ref{cnk}) we write for the coefficients $g_{\alpha,\beta}$ the recurrence formula
	\[g_{\alpha,\beta}=\frac{\alpha+\beta-1}{\alpha+\beta}g_{\alpha-1,\beta}+\frac{s-\beta+1}{\alpha+\beta}g_{\alpha-1,\beta-1}.
\]
Thus formula (\ref{Gam}) can be rewritten in the following form
\begin{equation}\label{Gam2}
	\Gamma (s+1)=\frac{1}{s+1}+\sum\limits_{\alpha=1}^{\infty}\frac{1}{\alpha +s+1}\sum\limits_{\beta=1}^{\alpha}g_{\alpha,\beta},
\end{equation}
which is more suitable for numerical computations. The convergence of the series (\ref{Gam}) and (\ref{Gam2}) will be discussed in section \ref{ne}.
\section{Derivative polynomials}
Let $\{a_1,a_2,\ldots ,a_n\}$ be a permutation of the set $\{1,2,\ldots ,n\}$. Then $\{a_{j},a_{j+1}\}$ is an ascent of the permutation if $a_j< a_{j+1}$. 
The Eulerian number $\displaystyle {n \bangle k} $ is defined as the number of permutations of the set $\{1,2,\ldots ,n\}$ having $k$ permutation ascents (see \cite{GKP}, p.267). For example for $n=3$ the permutation $\{1,2,3\}$ has two ascents, namely $\{1,2\}$ and $\{2,3\}$, and  $\{3,2,1\}$ has no ascents. Each of the other four permutations of the set has exactly one ascent. Thus $\displaystyle {3 \bangle 0} =1 $,  $\displaystyle {3 \bangle 1} =4 $, and $\displaystyle {3 \bangle 2} =1 $.\\
Consider a function $x=x(t)$ which satisfies Riccati's differential equation with constant coefficients 
\begin{equation}\label{Ric}
	x'(t)=ax^{2}+bx+c = a(x-\alpha)(x-\beta),
\end{equation}
where $a,b,c$ are real numbers, $a\neq 0$ and the roots $\alpha, \beta $ are real or complex conjugate numbers. Examples of such functions and equations are:
\begin{enumerate}
	\item $x(t)=\tan t, \quad x'(t)=x^{2}+1$,
	\item $x(t)=\tanh t, \quad x'(t)=-x^{2}+1$,
	\item $x(t)=1/(1+e^{ t}), \quad x'(t)=x^{2}-x$.
	\end{enumerate}
The following formula has been discussed during the Conference ICNAAM 2006 (September 2006) in Greece and it appears, with an inductive proof, in the paper \cite{Rz}. Independently the formula has been considered and proved, with a proof based on generating functions, by Franssens \cite{F} (see also \cite{R1}).\\
If a function $x(t)$ satisfies equation (\ref{Ric}), 
then the $n$th derivative of $x(t)$ can be expressed by the following formula:
\begin{equation}\label{f}
	x^{(n)}(t)=a^{n}\sum\limits_{k=0}^{n-1}{n \bangle k} 
(x-\alpha)^{k+1}(x-\beta)^{n-k},
\end{equation}
where $n=2,3,\ldots $.\\
The polynomials of the variable $x$ on the right hand side of (\ref{f}) are a particular case of the derivative polynomials which were introduced by Hoffman \cite{Ho} and recently intensively studied (see for example \cite{Bo}, \cite{Bo2}, \cite{F},  \cite{R2}, \cite{R3}).
Let us denote by $Q_{n+1}(x)$ the derivative polynomial of the $n$th derivative of the function $x(t)=1/(1+e^{ t})$ (where $\alpha =0$, $\beta =1$), of degree $n+1$. Thus for $n\ge 2$ 
	\[Q_{n}(x)=\sum\limits_{k=0}^{n-2}{n-1 \bangle k} 
x^{k+1}(x-1)^{n-1-k}.
\]
The polynomial $Q_{n}(x)$  ($n\ge 2$) is divisible by $x(x-1)$. We will denote by $P_{n-2}(x)$ the polynomial of degree $n-2$ resulting from the operation. Thus for any $n\ge 0$ we obtain
	\[P_{n}(x)=\sum\limits_{k=0}^{n}{n+1 \bangle k} 
x^{k}(x-1)^{n-k}.
\]
In the paper \cite{R1} (see also \cite{Bo}, \cite{Bo2}, \cite{F}) it is proved that the polynomials  $Q_{n}(x)$ can be expressed in terms of the Stirling numbers of 
the second kind ${ n \brace k}$ (number of the ways of partitioning a set of $n$ 
elements into $k$ nonempty subsets, see Graham et al. \cite{GKP}), namely
	\[Q_{n}(x)=\sum\limits_{k=1}^{n} (-1)^{n-k}(k-1)! {n \brace k} x^{k}.
\]
Then one can easily check that
\[P_{n}(x)=\sum\limits_{k=1}^{n+1} (-1)^{n+1-k}k! {n+1 \brace k} x^{k-1}.
\]
It will be useful for us to have the expansion of $P_{n}(x)$ at the point $x_0=1/2$. By rewriting the polynomial in the form with unknown coefficients 
	\[P_{n}(x)=\sum\limits_{j=0}^{n}p_{n,j}\left(x-\frac{1}{2}\right)^{j},
\]
one can verify that the following recurrence formula holds
	\[p_{n+1,j}=(j+1)\left(p_{n,j-1}-\frac{1}{4} p_{n,j+1}\right).
\]
From results of \cite{R1} it follows that for 
$n=1,2,\ldots $ the polynomial $P_n(x)$  has exactly $n$ 
simple zeroes which all are real and lie in the interval $(0,1)$. Moreover, the zeroes of the polynomials have the interlacing property. 
\section{The Riemann Zeta function}
It is well known that if $\textrm{Re}\: s>0$ then
\begin{equation}\label{Z}
	\int\limits_{0}^{\infty}t^{s-1}\frac{1}{1+e^{t}}dt=\zeta (s)(1-2^{1-s})\Gamma (s).
\end{equation}
Integrating $n+1$ times by parts the left hand side of (\ref{Z}) we obtain successively
\begin{eqnarray}
\lefteqn{\hspace{-15mm}\int\limits_{0}^{\infty}t^{s-1}\frac{1}{1+e^{t}}dt=
-\frac{1}{s}\int\limits_{0}^{\infty}t^{s}\left(\frac{1}{1+e^{t}}\right)'dt=
	\frac{1}{s(s+1)}\int\limits_{0}^{\infty}t^{s+1}\left(\frac{1}{1+e^{t}}\right)''dt
}\nonumber \\
&&=\cdots = \frac{(-1)^{n+1}}{s(s+1)\cdots (s+n)}\int\limits_{0}^{\infty}t^{s+n}\left(\frac{1}{1+e^{t}}\right)^{(n+1)}dt.\label{Zp}
\end{eqnarray}
Substituting in (\ref{Zp}) $x=1/(1+\exp(t))$ we obtain
{\setlength\arraycolsep{2pt}
\begin{eqnarray}
\zeta (s)(1-2^{1-s})\Gamma (s)&=& \frac{(-1)^{n}}{s(s+1)\cdots (s+n)}\int\limits_{0}^{1/2}\left( \log \frac{1-x}{x}\right)^{s+n}\frac{Q_{n+2}(x)}{x(x-1)}\;dx\nonumber \\
&=& \frac{(-1)^{n}}{s(s+1)\cdots (s+n)}\int\limits_{0}^{1/2}\left( \log \frac{1-x}{x}\right)^{s+n}P_{n}(x)\;dx,\label{Pn}
\end{eqnarray}}
where the polynomials $Q_{n}(x)$ and $P_{n}(x)$ have been defined in the previous section. \\
Our next aim is to find a power series expansion of the function $(\log((1-x)/x))^{r}$ by using the potential polynomials (\ref{pp2}). Since
	\[\log\frac{1-x}{x}=\log\frac{1+(1-2x)}{1-(1-2x)}=
	2(1-2x)\left( 1+\frac{1}{3}(1-2x)^{2} +\frac{1}{5}(1-2x)^{4}+\cdots \right),
\]
we have to find the expansion of the function
	\[\left(\frac{1}{2t}\log\frac{1+t}{1-t}\right)^{r}=\left( 1+\frac{1}{3}t^{2} +\frac{1}{5}t^{4}+\cdots \right)^{r}.
\]
In order to do this we will compute the Bell exponential polynomials\\ $B_{n,k}(0,2!/3,0,4!/5,0, 6!/7,\ldots)$. By using (\ref{Bp}) we get
\begin{eqnarray*}
\lefteqn{\hspace{-20mm}\exp\biggl\{u\left(\frac{1}{3}t^{2}+\frac{1}{5}t^{4}+\frac{1}{7}t^{6}+\cdots\right)\biggr\}
=e^{-u}\exp\biggl\{u\left(1+\frac{1}{3}t^{2}+\frac{1}{5}t^{4}+\frac{1}{7}t^{6}+\cdots\right)\biggr\}}
\\ && =e^{-u}\exp\biggl\{\frac{u}{2t}2\left(t+\frac{1}{3}t^{3}+\frac{1}{5}t^{5}+\frac{1}{7}t^{7}+\cdots\right)\biggr\}
\\ && =e^{-u}\exp\biggl\{\frac{u}{2t}\log\frac{1+t}{1-t}\biggr\}=e^{-u}\left(\frac{1+t}{1-t}\right)^{\frac{u}{2t}}.
\end{eqnarray*}
It is well known that the function $((1+t)/(1-t))^{x}$ is a generating function for the Mittag-Leffler polynomials $M_{k}(x)$
	\[\left(\frac{1+t}{1-t}\right)^{x}=\sum\limits_{k=0}^{\infty}\frac{M_{k}(x)}{k!}t^{k},\qquad |t|<1.
\]
A first few Mittag-Leffler polynomials are as follows
{\setlength\arraycolsep{2pt}
\begin{eqnarray*}
M_0(x)&=& 1,\;M_{1}(x)=2x,
\;M_{2}(x)=4x^{2},\; M_{3}(x)=8x^{3}+4x,\\
M_{4}(x)&=&16x^{4}+32x^{2},\;M_{5}(x)=32x^{5}+160x^{3}+48x.
\end{eqnarray*}}
Bateman \cite{Ba} considers polynomials $g_{k}(x)=M_{k}(x)/k!$ and gives for them the following recurrence formula
\begin{equation}\label{gn}
	ng_{n}(x) = (n - 2)g_{n-2}(x) + 2xg_{n-1}(x),
\end{equation}
with which he refers to Belorizky \cite{Be}. Multiplying both sides of (\ref{gn}) by $(n-1)!$ we get
\begin{equation}\label{Mn}
M_{n}(x)=(n-1)(n-2)M_{n-2}(x)+2xM_{n-1}(x).
\end{equation}
If  $M_{n}(x)=\sum\limits_{k=0}^{n}a_{n,k}x^{k} $ then the recurrence formula (\ref{Mn}) yields
\begin{equation}\label{ank}
	a_{n,k}=(n-1)(n-2)a_{n-2,k}+2a_{n-1,k-1}.
\end{equation}
\begin{lemma}
Numbers $a_{n,k}$ fulfill the following formula
\begin{equation}\label{anka}
	a_{n,k}=\sum\limits_{m=0}^{k-1}2^{m}(n-1-m)(n-2-m)a_{n-2-m,k-m}.
\end{equation}
\end{lemma}
\begin{proof}
Formula (\ref{anka}) follows by using formula (\ref{ank}) in each step to the last term of the resulting sum.
\end{proof}
Therefore
\begin{eqnarray*}
\lefteqn{e^{-u}\left(\frac{1+t}{1-t}\right)^{\frac{u}{2t}}=
e^{-u}\sum\limits_{n=0}^{\infty}\frac{M_{n}(u/2t)}{n!}t^{n}
=e^{-u}\sum\limits_{n=0}^{\infty}\frac{1}{n!}
\sum\limits_{k=0}^{n}a_{n,k}\left(\frac{u}{2t}\right)^{k}t^{n}
}
\\ && =\left(\sum\limits_{j=0}^{\infty}(-1)^{j}\frac{u^{j}}{j!}\right)\sum\limits_{n=0}^{\infty}
\sum\limits_{k=0}^{n}\frac{1}{n!}\;a_{n,k}\;\frac{u^{k}}{2^{k}}\;t^{n-k}
\\&& =\sum\limits_{\alpha=0}^{\infty}\sum\limits_{\beta=0}^{\alpha}\sum\limits_{j=0}^{\beta}
(-1)^{j}\frac{1}{j!}\frac{1}{(\alpha+\beta-j)!}\frac{1}{2^{\beta-j}}\;a_{\alpha+\beta-j,\beta-j}\;u^{\beta}t^{\alpha}
\\&& =\sum\limits_{\alpha=0}^{\infty}\sum\limits_{\beta=0}^{\alpha}\frac{1}{(\alpha+\beta)!}\sum\limits_{j=0}^{\beta}
(-1)^{j}{\alpha+\beta\choose j}\frac{1}{2^{\beta-j}}\;a_{\alpha+\beta-j,\beta-j}\;u^{\beta}t^{\alpha},
\end{eqnarray*}
and we see that, in our case, the value of the Bell polynomial $B_{\alpha,\beta}$   is
\begin{equation}\label{Bell}
	B_{\alpha,\beta}=\frac{\alpha!}{(\alpha+\beta)!}\sum\limits_{j=0}^{\beta}
	(-1)^{j}{\alpha+\beta\choose j}\frac{1}{2^{\beta-j}}\;a_{\alpha+\beta-j,\beta-j}.
\end{equation}
Let us denote
\begin{equation}\label{b}
	b_{\alpha,\beta}=\sum\limits_{j=0}^{\beta}
	(-1)^{j}{\alpha+\beta\choose j}\frac{1}{2^{\beta-j}}\;a_{\alpha+\beta-j,\beta-j}.
\end{equation}
The following Table \ref{tab:tb2} gives the first few values of $b_{\alpha,\beta}$.

\begin{table}[ht]
\caption{Coefficients $b_{\alpha,\beta}$}
\centering
        $\begin{array}[]{ccccccc}
        \alpha\backslash\beta &\quad 0 & \quad 1 & 2 & 3 & 4 & \quad 5   \\
        0 &\quad 1 &\quad 0 & &&&\\
        1 &\quad 0 &\quad 0 &  &&& \\
        2 &\quad 0 &\quad 2 & 0 &&& \\
        3 &\quad 0 &\quad 0 & 0 && &\\
        4 &\quad 0 &\quad 24 & 40 & 0 & &\\
        5 &\quad 0 &\quad 0 & 0 & 0 & & \\
        6 &\quad 0 &\quad 720 & 2688 & 2240 & 0 &  \\
        7 &\quad 0 &\quad 0 & 0 & 0 & 0 &  \\
        8 &\quad 0 &\quad 40320 & 245376 & 443520 & 246400 & \quad 0.
        \end{array}
$
\label{tab:tb2}
\end{table}
\begin{lemma}
Numbers $b_{\alpha,\beta}$ fulfill the recurrence formula 
\begin{equation}\label{bnk}
b_{\alpha,\beta}=(\alpha+\beta-2)(\alpha+\beta-1)(b_{\alpha-2,\beta}+b_{\alpha-2,\beta-1}).
\end{equation}
\end{lemma}
\begin{proof}
Assume that $\alpha\ge 4$. Since the last term in (\ref{b}) is zero then as the upper limit of the sum we can take $\beta-1$.  We have
{\setlength\arraycolsep{2pt}
\begin{eqnarray}
 b_{\alpha-2,\beta}&=&\sum\limits_{j=0}^{\beta-1}
	(-1)^{j}{\alpha+\beta-2\choose j}\frac{1}{2^{\beta-j}}\;a_{\alpha+\beta-2-j,\beta-j},\label{b1}\\
 b_{\alpha-2,\beta-1}&=&\sum\limits_{j=0}^{\beta-2}
	(-1)^{j}{\alpha+\beta-3\choose j}\frac{1}{2^{\beta-1-j}}\;a_{\alpha+\beta-3-j,\beta-j}\label{b2}.
\end{eqnarray}}
By adding similar terms in (\ref{b1}) and (\ref{b2}) we see that formula (\ref{bnk}) is equivalent to 
\begin{equation}\label{bnk1}
	b_{\alpha,\beta}=(\alpha+\beta-2)(\alpha+\beta-1)\sum\limits_{j=0}^{\beta-1}
	(-1)^{j}{\alpha+\beta-3\choose j}\frac{1}{2^{\beta-j}}\;a_{\alpha+\beta-2-j,\beta-j}.
\end{equation}
From the other side
{\setlength\arraycolsep{2pt}
\begin{eqnarray*}
b_{\alpha,\beta}&=&\sum\limits_{j=0}^{\beta-1}
	(-1)^{j}{\alpha+\beta\choose j}\frac{1}{2^{\beta-j}}\;a_{\alpha+\beta-j,\beta-j}=
	\sum\limits_{j=0}^{\beta-1} (-1)^{j}{\alpha+\beta\choose j}\frac{1}{2^{\beta-j}}\times\\
	&&\times\sum\limits_{i=0}^{\beta-j-1}2^{i}(\alpha+\beta-j-i-1)(\alpha+\beta-j-i-2)a_{\alpha+\beta-j-i-2,\beta-j-i}\\
	&&\hspace{-10mm}=\sum\limits_{m=0}^{\beta-1}\frac{1}{2^{\beta-m}}
	(\alpha+\beta-m-1)(\alpha+\beta-m-2)a_{\alpha+\beta-m-2,\beta-m}
	\sum\limits_{k=0}^{m}(-1)^{k}{\alpha+\beta\choose k}\\
	&&\hspace{-10mm}=\sum\limits_{m=0}^{\beta-1}\!\frac{(-1)^{m}}{2^{\beta-m}}{\alpha\!+\!\beta\!-\!1\choose m}
	(\alpha\!+\!\beta\!-\!m\!-\!1)(\alpha\!+\!\beta\!-\!m\!-\!2)a_{\alpha+\beta-m-2,\beta-m}\\
	&&\hspace{-10mm}=(\alpha\!+\!\beta\!-\!1)(\alpha\!+\!\beta\!-\!2)\sum\limits_{m=0}^{\beta-1}\frac{(-1)^{m}}{2^{\beta-m}}
	{\alpha\!+\!\beta\!-\!3\choose m}a_{\alpha+\beta-m-2,\beta-m},
\end{eqnarray*}}
where $m=i+j$. Thus formula (\ref{bnk}) is proved. At the beginning of the above calculation we used formula (\ref{anka}) and then the following formulae 
{\setlength\arraycolsep{2pt}
\begin{eqnarray*}
&&\sum\limits_{k=0}^{m}(-1)^{k}{n \choose k}=(-1)^{m}{n-1 \choose m},\\
&&{n-1\choose m}(n-m-1)(n-m-2)={n-3\choose m}(n-1)(n-2),
\end{eqnarray*}}
which are easy to check.
\end{proof}
Using formula (\ref{pp}) and the values of the Bell polynomials (\ref{Bell}) we get
\begin{eqnarray}
\lefteqn{\left(\log\frac{1-x}{x} \right)^{r}=\left(\log\frac{1+(1-2x)}{1-(1-2x)} \right)^{r}
}\nonumber\\
&&=2^{r}(1-2x)^{r}\left(1+\frac{(1-2x)^{2}}{3}+\frac{(1-2x)^{4}}{5}+\cdots\right)^{r}\nonumber\\
&&=2^{r}(1-2x)^{r}\left\{1+\frac{2}{3!}(1-2x)^{2}r+(1-2x)^{4}\left(\frac{24}{5!}r+\frac{40}{6!}r(r-1)\right)\right.\nonumber\\
&&\left. +(1-2x)^{6}\left(\frac{720}{7!}r+\frac{2688}{8!}r(r-1)+\frac{2240}{9!}r(r-1)(r-2)\right)+\cdots
  \right\}\nonumber\\
&&=2^{r}(1-2x)^{r}\biggl(1+\sum\limits_{m=1}^{\infty}(1-2x)^{2m}\sum\limits_{k=1}^{m}\frac{b_{2m,k}}{(2m+k)!}(r)_{k}\biggr). \label{G1} 
\end{eqnarray}
Therefore since
	\[\int\limits_{0}^{1/2}(1-2x)^{\gamma}dx=\frac{1}{2(\gamma+1)},
\]
then, putting in (\ref{Pn}) say $n=0$ ($P_{0}(x)=1$) we get
\begin{eqnarray}
\lefteqn{\hspace{-5mm}\zeta (s)(1-2^{1-s})\Gamma (s)= \frac{2^{s-1}}{s}\left\{\frac{1}{s+1}+\frac{2s}{3!(s+3)}+\frac{1}{s+5}\left(\frac{24}{5!}s+\frac{40}{6!}s(s-1)\right)\right.
}\nonumber\\
&&\left. +\frac{1}{s+7}\left(\frac{720}{7!}s+\frac{2688}{8!}s(s-1)+\frac{2240}{9!}s(s-1)(s-2)\right)+\cdots
  \right\}\nonumber\\
&&=\frac{2^{s-1}}{s}\biggl(\frac{1}{s+1}+\sum\limits_{m=1}^{\infty}\frac{1}{s+2m+1}\sum\limits_{k=1}^{m}\frac{b_{2m,k}}{(2m+k)!}(s)_{k}\biggr).\label{G2} 
\end{eqnarray}
By denoting here
	\[z_{m,k}=\frac{b_{m,k}}{(m+k)!}(s)_{k},
\]
and using the recurrence formula for $b_{m,k}$ we see that coefficients $z_{m,k}$ fulfill the recurrence formula
	\[z_{m,k}=\frac{m+k-2}{m+k}z_{m-2,k}+\frac{s-k+1}{m+k}z_{m-2,k-1}.
\]
Therefore the above expansion can be rewritten in the form
\begin{equation}\label{G3}
	\zeta (s)(1-2^{1-s})\Gamma (s)=\frac{2^{s-1}}{s}\biggl(\frac{1}{s+1}+\sum\limits_{m=1}^{\infty}\frac{1}{s+2m+1}\sum\limits_{k=1}^{m}z_{2m,k}\biggr),
\end{equation}
which is more useful for numerical computations. The convergence of the series (\ref{G2}) (and (\ref{G3})) is investigated in section \ref{ne}.
\section{Numerical experiments}\label{ne}
\begin{figure}[!ht]
        \centering
       \includegraphics{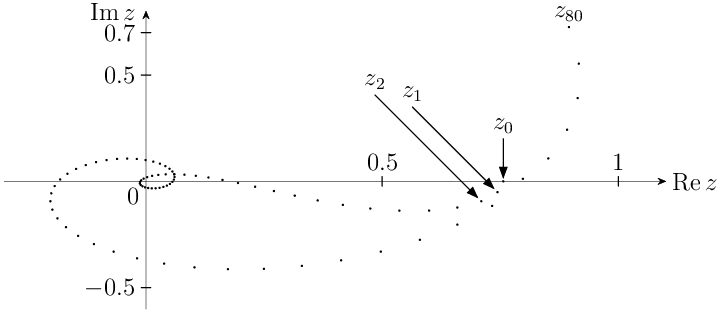}
        \caption{The partial sum of the series (\ref{G4}) at points $s_{k}=1/2+ik/5; \; k=0,1,\ldots, 80$. }
        \label{fig:z2}
\end{figure}

Flajolet and Odlyzko (\cite{FO}, Theorem 3A page 227) proved that if $a$, $b$ are any complex numbers and $a,b \notin \{0,1,2,\ldots \}$ then the Taylor coefficients $\{f_{n}\}$ of the function (analytic on the open unit disk and having one singularity at $z=1$ on the unit circle)
	\[f(z)=(1-z)^{a}\left(
	\frac{1}{z}\log \frac{1}{1-z}\right)^{b}
\]
satisfy
\begin{equation}\label{ex1}
	f_{n}=[z^{n}]f(z)\sim \frac{n^{-a-1}}{\Gamma (-a)}(\log n)^{b}\left(1+\sum\limits_{k\ge 1}\frac{e_{k}^{(a,b)}}{\log^{k} n}\right)
\end{equation}
with
	\[e_{k}^{(a,b)}=(-1)^{k}{b \choose k } \Gamma(-a)\frac{d^{k}}{ds^{k}}\left(\frac{1}{\Gamma (-s)}\right)\!\Biggl|_{s=a}\biggr..
\]
How to treat the limit case when, for example, $a$ is a nonnegative integer is shown briefly later in this paper (see the remarks on pp. 231--232). There is a slightly expanded but entirely similar discussion 
in the book \cite{FS} by Flajolet and Sedgewick (see Theorem VI.2 page 385 and the remarks on pp. 386--387). From this it follows, for instance, that expression (\ref{ex1}) has a limit when $a$ goes to zero (the first term is zero because we take $1/\Gamma(0)=0$ and in the next terms the quantity $\Gamma(-a)$ is cancelled). Therefore we get the following estimate for the $n$th Taylor coefficient of the function $(-\log(1-z)/z)^{b}$
\begin{equation}\label{es}
	[z^{n}]\left(
	\frac{1}{z}\log \frac{1}{1-z}\right)^{b} = O\left(\frac{b(\log n)^{b-1}}{n} \right).
\end{equation}
The coefficient of $t^{\alpha +s} $ in (\ref{Ga}) 
	\[\sum\limits_{\beta=1}^{\alpha}(s)_{\beta}\frac{c_{\alpha,\beta}}{(\alpha+\beta)!}
\]
is the $\alpha$th Taylor coefficient of the function $(-\log(1-t)/t)^{s}$. By (\ref{es}) we see that the series (\ref{Gam}) (and therefore (\ref{Gam2})) are absolutely convergent and the terms of the series have the estimate
\begin{equation}
	\frac{1}{\alpha +s+1}\sum\limits_{\beta=1}^{\alpha}(s)_{\beta}\frac{c_{\alpha,\beta}}{(\alpha+\beta)!}=
	O\left(\frac{|s|(\log \alpha)^{\textrm{\scriptsize{Re}} s-1}}{\alpha \;|\alpha +1+s|} \right),
\end{equation}
 as $\alpha \rightarrow \infty$.
Similar reasoning can be made for the function 
	\[f(z)=\left(
	\frac{1}{2z}\log \frac{1+z}{1-z}\right)^{b},
\]
which is analytic on the open unit disk and has two singularities (at $z=1$ and $z=-1$) on the unit circle. Flajolet and Sedgewick (\cite{FS}, Theorem VI.5 page 398) also discuss the case of multiple singularities and prove that the contributions from each of the singularities are to be added. From this it follows that the rate of convergence of the series (\ref{G2}) (and (\ref{G3})) is the same as in the previous case.\\
We will show now the results of two numerical experiments. The first is associated with formula (\ref{Gam2}). The Table \ref{tab:tb3} below gives the values of the partial sums of  (\ref{Gam2}) for $s=-1/2$. The exact value of the sum of the series is $\Gamma(1/2)=\sqrt{\pi}=1.77245\ldots $. 
\begin{table}[h]
\caption{The partial sums of  (\ref{Gam2}) for $s=-1/2$}
\centering
$\begin{array}[]{|c|c|}
\hline
         n= & n\textrm{th partial sum of}\;(\ref{Gam2})  \\
         \hline \hline
        20 & 1.77460\ldots  \\
        \hline
        50 & 1.77315\ldots \\
        \hline
        73 &1.77289\ldots \\    
        \hline
        \end{array}
$
\label{tab:tb3}
\end{table}
These calculations confirm the, above obtained, estimates of the rate of convergence of the expansion (\ref{Gam2}).\\
The second experiment is related to formula (\ref{G3}). We have computed the values of the 35th (the upper limit of summation for $m$ is 35)  partial sum of the series
\begin{equation}\label{G4}
\frac{1}{s+1}+\sum\limits_{m=1}^{\infty}\frac{1}{s+2m+1}\sum\limits_{k=1}^{m}z_{2m,k}=
s(1-2^{1-s})2^{1-s}\Gamma (s)\zeta (s)
\end{equation}
at 81 equally spaced points of the interval $\{s=1/2+it\}; \;\; 0\le t \le 16$. The result has been shown on Figure \ref{fig:z2}, where $z_{k}$ is the value of the partial sum at the point $s=s_{k}=1/2+ik/5; \; k=0,1,\ldots ,80$.
\section{Concluding remarks}
Expansions similar to (\ref{G3}) can be established for any $n\ge 1$ in (\ref{Pn}). In fact the expansion can be seen as an expansion for the integral
\begin{equation}\label{int}
	2^{1-s-n}(-1)^{n}\int\limits_{0}^{1/2}\left( \log \frac{1-x}{x}\right)^{s+n}\!\!\!P_{n}(x)\;dx=\frac{s^{(n+1)} }{2^{s+n-1}}\;\zeta (s)(1-2^{1-s})\Gamma (s),
\end{equation}
where $s^{(n+1)}=s(s+1)\cdots (s+n)$ denotes the rising factorial.\\
If we choose $s\in C$ such that $1/2<\textrm{Re}\: s\le 1$ and $\zeta (s)\neq 0$ then it can be seen from (\ref{int}) that the left hand side (i.e. our expansion) goes to infinity for $n\rightarrow \infty$. Therefore if for any $s\in C$ ($1/2 < \textrm{Re}\: s < 1$) we had indicated a nonnegative integer $n$ such that the evaluation, at this point $s$, of the expansion of the left hand side of (\ref{int}) is different from zero then the RH would follow.
\section{Acknowledgments}
I would like to thank P. Flajolet for pointing out to me the formulae used in Section 5, and for comments on them. I am also grateful to M. Skwarczy\'nski for many valuable discussions during the work on this paper.

\end{document}